\documentclass[11pt]{amsart}
\usepackage{amsmath}
\usepackage{amsfonts}
\usepackage{amssymb}
\newtheorem{thm}{Theorem}[section]

\newtheorem{lemma}[thm]{Lemma}

\theoremstyle{definition}

\newtheorem{remark}[thm]{Remark}

\newcommand{\bb}[1]{\mathbb{#1}}
\newcommand{\cl}[1]{\mathcal{#1}}

\begin{document}

\title[Syndetic Sets]{Syndetic Sets, Paving and the Feichtinger Conjecture}

\author[V.~I.~Paulsen]{Vern I.~Paulsen}
\address{Department of Mathematics, University of Houston,
Houston, Texas 77204-3476, U.S.A.}
\email{vern@math.uh.edu}

\thanks{This research was supported in part by NSF grant DMS-0600191 and by
  the American Institute of Mathematics.}

\begin{abstract}
We prove that if a Bessel
sequence in a Hilbert space, that is indexed by a countably infinite
group in an invariant manner, can be partitioned
into finitely many Riesz basic sequences, then each of the sets in the
partition can be chosen to be syndetic.  We then apply this result to
prove that if a Fourier frame for a measurable subset of a higher
dimensional cube can be partitioned into Riesz basic sequences, then
each subset can be chosen to be a syndetic subset of the corresponding
higher dimensional integer lattice. Both of these results follow from a result about syndetic pavings of elements of the von~Neumann algebra of a discrete group.
\end{abstract}

\maketitle


\section{Introduction} 

Recently, W. Lawton\cite{La} proved that if $A \subseteq [0,1]$ is a
set of positive Lebesgue measure such that the Fourier frame $\{ e^{2 \pi i
  nt}\chi_A \}_{n \in \bb Z}$ for $L^2(A)$ can be partitioned into
Riesz basic sequences, then each set in the partition can be chosen
to be a syndetic subset of $\bb Z.$ In my dynamical systems approach
to the Kadison-Singer problem\cite{Pa}, I proved that an operator in
the von~Neumann algebra of a countable abelian group that could
be paved in the sense of Anderson, could be paved by syndetic
sets. Although this fact was derived in the proof of
\cite[Theorem~16]{Pa}, it is not apparent from the statement of the
theorem and has only been announced without proof at various meetings. If one combines this result about syndetic paving of elements of the
von~Neumann algebra of a group with the method of proof of the equivalence
of the Kadison-Singer conjecture and the Feichtinger conjecture\cite{CT}, then
one arrives, in this roundabout way, at a result that generalizes
Lawton's result.

In this paper we make this argument direct. We begin with a direct
proof of the result about syndetic paving of elements of von~Neumann
algebras of groups(Theorem~\ref{pavingthm}).  This direct proof has the added
advantage that we are able to eliminate the
hypothesis that the group is abelian. We then show how this paving
result directly yields statements about partitioning certain Bessel sequences
into Riesz basic sequences. In particular, in addition to the paving
theorem, we shall prove
the following two theorems:

\begin{thm}\label{thm1} Let $G$ be a countable group, let $\cl H$ be a
  Hilbert space and let $\{ f_g \}_{g \in G}$ be a Bessel sequence in
  $\cl H$ with the property that for any $g_1,g_2,h_1,h_2 \in G,$
  satisfying $g_1g_2^{-1} = h_1h_2^{-1}$ we
  have $\langle f_{g_1}, f_{g_2} \rangle = \langle f_{h_1},f_{h_2}
  \rangle .$  If $G$ can be partitioned into finitely many subsets $\{
  A_l: 1 \le l \le L \}$ such that for $1 \le l \le L,$ $\{ f_g: g \in A_l \}$ is a Riesz
  basic sequence, then $G$ can be partitioned into $K$ syndetic
  subsets $\{ S_k: 1 \le k \le K \}$ with $K \le L$ such that for $1
  \le k \le K,$ $\{ f_g:
  g \in S_k \}$ is a Riesz basic sequence.
\end{thm}

To state the second theorem it wil be convenient to introduce
some notation.  Given a natural number $N,$ we let $[0,1]^N$
denote the product of $N$ copies of the unit interval and we let $\bb
Z^N$ denote the group that is the direct sum of $N$ copies of the
group of integers. If $t=(t_1,...,t_N) \in [0,1]^N$ and
$n=(n_1,...,n_N) \in \bb Z^N,$ then we set $n \cdot t = n_1t_1+ \cdots
+n_Nt_N$ and set $e_n(t) = e^{2 \pi i n \cdot t }.$ 
We have that $\{ e_n(t): n \in \bb Z^N \}$ is an orthonormal basis for
the Hilbert space $L^2([0,1]^N)$ of square-integrable functions with
respect to Lebesgue product measure. If $A \subseteq [0,1]^N$ is a set
of positive Lebesgue measure and $L^2(A)$ denotes the corresponding
Hilbert space of square-integrable functions on $A,$ then $f_n(t) =
e_n(t) \chi_A, n \in \bb Z^N$ is a Parseval frame for $L^2(A)$ called
the {\em Fourier frame for $L^2(A).$}

\begin{thm}\label{thm2} Let $N$ be a natural number and let $A
  \subseteq [0,1]^N$ be a set of positive Lebesgue measure. If the
  Fourier frame $\{ f_n \}_{n \in \b Z^N}$ for $L^2(A)$ can be
  partitioned into finitely many subsets $\{
  A_l: 1 \le l \le L \}$ such that for $1 \le l \le L,$ $\{ f_n: n \in A_l \}$ is a Riesz
  basic sequence, then $\bb Z^N$ can be partitioned into $K$ syndetic
  subsets $\{ S_k: 1 \le k \le K \}$ with $K \le L$ such that for $1
  \le k \le K,$ $\{ f_n:
  n \in S_k \}$ is a Riesz basic sequence.  
\end{thm}

\section{Proofs of the Main Results}

We begin by recalling some results and notation from \cite{Pa}.
Let $G$ be a countably infinite group, let $\ell^2(G)$ denote the Hilbert space of square-summable functions on $G,$ let $\ell^{\infty}(G)$ denote the algebra of bounded functions on $G,$ let $\beta G$ denote the Stone-Cech compactification of $G$ and let $C(\beta G)$ denote the algebra of continuous functions on $\beta G.$ Every element of $\ell^{\infty}(G)$ extends uniquely to a continuous function on $\beta G,$ and conversely, and hence, we identify $\ell^{\infty}(G) = C(\beta G).$  We let $e_g \in \ell^2(G)$ denote the function that is $1$ at $g$ and $0$ elsewhere, so that $\{ e_g: g \in G \}$ is an orthonormal basis for $\ell^2(G).$ We let $B(\ell^2(G))$ denote the space of bounded operators on $\ell^2(G)$ and given a bounded operator $X$ we identify $X$ with a $G \times G$ matrix by setting $X= (x_{g,h})$ where $x_{g,h} = \langle Xe_h, e_g \rangle .$ If a bounded operator is diagonal with respect to this basis, then its diagonal entries define an element of $\ell^{\infty}(G),$ so we also identify $\ell^{\infty}(G) \subseteq B(\ell^2(G))$ as the abelian subalgebra of diagonal operators. We let $E:B(\ell^2(G)) \to \ell^{\infty}(G)$ denote the completely positive projection defined by setting $E(X)$ equal to the diagonal part of the matrix $X.$

We let $\lambda:G \to B(\ell^2(G))$ denote the left regular representation of $G$ given by $\lambda(g)e_h = e_{gh}.$ The von~Neumann algebra of the group is the double commutant of the left regular representation, $VN(G) = \lambda(G)^{\prime \prime}.$ 

Given $X=(x_{g,h}) \in B(\ell^2(G)),$ we set $D_g = E(X\lambda(g^{-1})),$ that is, $D_g$ is the diagonal matrix $(y_{h,h})$ where $y_{h,h} = x_{h, g^{-1}h}.$ The operators $D_g$ can be thought of as the terms in a {\em formal series for $X,$} in the sense that $\sum_{g \in G} D_g\lambda(g)$ converges entrywise to the matrix for $X.$ We shall write $X \sim \sum_{g \in G} D_g \lambda(g).$ It can be easily seen that $X \in VN(G)$ if and only if $X \sim \sum_{g \in G} c_g \lambda(g)$ for some constants $c_g.$ That is, $X \in VN(G)$ if and only if $E(X\lambda(g^{-1}))$ is a scalar multiple of the identity operator for every $g \in G.$

We also need to recall the semigroup structure of $\beta G.$ For a good reference, see the text of Hindman-Strauss\cite{HS}, but we also recall the salient points below.
The action of $G$ on $G$ given by left translation extends to give an action of $G$ on $\beta G.$ In particular, if $g \in G,$ $\omega \in \beta G,$ and $\{ g_{\lambda} \}$ is a net in $G$ converging to $\omega,$ then the net $\{ gg_{\lambda} \}$ converges to the element $g \cdot \omega.$  Moreover, this left action extends to define an associative product making $\beta G$ into a right continuous, compact semigroup. If $\alpha, \omega \in \beta G,$ and $\{ g_{\lambda} \}, \{ h_{\mu} \}$ are nets in $G$ with $\alpha = \lim_{\mu} h_{\mu},$ $\omega = \lim_{\lambda} g_{\lambda},$ then $\alpha \cdot \omega = \lim_{\mu}[ \lim_{\lambda}   h_{\mu}g_{\lambda}].$

Fixing an element $\omega \in \beta G,$ we defined in \cite{Pa} a *-homomorphism, $\pi_{\omega}: C(\beta G) \to C(\beta G),$ by setting $\pi_{\omega}(f)(g) = f(g \cdot \omega).$ Via our identifications, we also write $\pi_{\omega}(D),$ where $D \in \ell^{\infty}(G) \subseteq B(\ell^2(G))$ is the diagonal operator that is identified with the function $f.$ In \cite{Pa}, we proved that there is a well-defined, unital completely positive map $\psi_{\omega}: B(\ell^2(G)) \to B(\ell^2(G))$ defined as follows:
\[\text{if } X \sim \sum_{g \in G} D_g \lambda(g), \text{ then } \psi_{\omega}(X) \sim \sum_{g \in G} \pi_{\omega}(D_g)\lambda(g). \]

Given a subset $A \subseteq G,$ we let $P_A$ denote the diagonal operator $(y_{h,h})$ with $y_{h,h}=1$ when $h \in A$ and $0,$ otherwise. That is, $P_A$ is the diagonal operator that is identified with the characteristic function of the set $A,$ $\chi_A.$

In \cite[Theorem~2.8]{HS} it is shown that the compact right topological semigroup, $\beta G$ contains a non-empty minimal(two sided) ideal, denoted $K(\beta G).$ Finaly, recall that a subset $S \subseteq G$ is called {\bf syndetic} if there exists a finite set $F \subseteq G,$ such that $G= \cup_{g \in F} g^{-1}S,$ i.e., $G = F^{-1}\cdot S.$

\begin{lemma}  Let $G$ be a countably, infinite discrete group, let $p \in K(\beta G),$ and let $A \subseteq G.$ If $p \in K(\beta G),$ then $\pi_p(P_A) = P_R,$ with $R$ either empty or syndetic.
\end{lemma}
\begin{proof} Since $\pi_p$ is a *-homorphism, $\pi_p(P_A)$ is a projection in $\ell^{\infty}(G)$ and hence, $\pi_p(P_A) = P_R$ for some subset $R \subseteq G.$
If $R$ is non-empty, then there exists $g_0 \in R$ and hence $1= \pi_p(P_A)(g_0) = f(g_0 \cdot p),$ where $f$ denotes the continuous of $P_A$ to $\beta G.$  That is, $f$ is the characteristic function of the closure of $A,$ $A^-$ in $\beta G.$ Thus, $g_0 \cdot p \in A^-,$ or $p \in (g_0^{-1}A)^-.$
Thus, we have that $R = \{ g \in G: p \in (g^{-1}A)^- \}.$

Using the identification of points in $\beta G$ with their corresponding ultrafilters, we see that this means that the set $g_0^{-1}A$ belongs to the ultrafilter $p,$ i.e., $g_0^{-1}A \in p.$ In this notation, $R = \{ g \in G: g^{-1}A \in p \}.$
  
Since $p \in K(\beta G)$ and $g_0^{-1}A \in p,$ by \cite[Theorem~4.39]{HS} $S= \{ x \in G: x^{-1}g_0^{-1}A \in p \}$ is syndetic.  Hence, $R = g_0S$ is also syndetic.
\end{proof}

We shall also need a consequence of a result of Choi's from the theory of multiplicative domains\cite{Ch}.

\begin{lemma}[Choi] Let $\cl A$ and $\cl B$ be unital C*-algebras and let $\phi: \cl A \to \cl B$ be a unital, completely positive map. If $1_{\cl A} \in \cl C \subseteq \cl A$ is a C*-subalgebra such that the restriction of $\phi$ to $\cl C$ is a *-hommomorphism $\pi,$ then for any $c_1,c_2 \in \cl C$ and any $a \in \cl A,$ we have that $\phi(c_1ac_2) = \pi(c_1) \phi(a) \pi(c_2).$
\end{lemma}
\begin{proof} One can regard $\cl B \subseteq B(\cl H)$ for some Hilbert space $\cl H$ so that Stinespring's representation theorem can be invoked. This implies the existence of a Hilbert space $\cl K,$ an isometry $V: \cl H \to \cl K$ and a *-homomorphism, $\rho: \cl A \to B(\cl K),$ such that $\phi(a) = V^*\rho(a)V.$ Restricting to $\cl C,$ one sees that $V \cl H$ reduces $\rho$ and that $\rho(c) = V\pi(c)V^* \oplus \pi_2(c)$ for some *-homomorphism $\pi_2.$ The result now follows.
\end{proof}

We can now prove a result about paving elements of the von~Neumann algebra of a group. For the case that $G$ is abelian, this result was announced at the 2008 Cincinnati GPOTS.

\begin{thm}\label{pavingthm} Let $G$ be a countably infinite discrete group, let $X \in VN(G)$ with $E(X) = 0,$ and let $\epsilon < 1.$  If $G$ can be partitioned into finitely many subsets $\{
  A_l: 1 \le l \le L \}$ such that for $1 \le l \le L,$ $\| P_{A_l}XP_{A_l}\| \le \epsilon \|X\|,$ then $G$ can be partitioned into $K$ syndetic
  subsets $\{ S_k: 1 \le k \le K \}$ with $K \le L$ such that for $1
  \le k \le K,$ $\| P_{S_k}XP_{S_k} \| \le \epsilon \|X\|.$ 
\end{thm}
\begin{proof}  Fix any $p \in K(\beta G),$ then we have a unital, completely positive map, $\psi_p: B(\ell^2(G)) \to B(\ell^2(G)).$ Since such maps are contractive, we have that $\|\psi_p(P_{A_l}XP_{A_l})\| \le \epsilon \|X\|.$  Note that if $D \in \ell^{\infty}(G)$ is a diagonal operator, then $\psi_p(D) = \pi_p(D),$ and so by Choi's Lemma, $\psi_p(P_{A_l}XP_{A_l}) = \pi_p(P_{A_l}) \psi_p(X) \pi_p(P_{A_l}).$

By the first lemma, $\pi_p(P_{A_l}) = P_{S_l},$ where $S_l$ is either
syndetic or empty.  Also, since $\pi_p$ is a unital *-homomorphism, $I
= \pi_p(P_{A_1} + \cdots + P_{A_L}) = P_{S_1} + \cdots + P_{S_L},$
which implies that $G$ is the disjoint union of the non-empty syndetic sets.

Finally, note that since $X \in VN(G),$ there exist scalars so that $X \sim \sum_{g \in G} c_gI \lambda(g),$ and hence, $\psi_p(X) \sim \sum_{g \in G} \pi_p(c_gI) \lambda(g) \sim X,$ that is, $\psi_p(X) = X.$
Thus, we have that $P_{S_l}XP_{S_l} = \psi_p(P_{A_l}XP_{A_l}),$ and the result follows.
\end{proof}

\begin{remark}
In the case that $G = \bb Z^N$ the Fourier transform defines an
isomorphism between $\ell^2(\bb Z^N)$ and $L^2([0,1]^N).$ That is,
choosing the orthonormal basis $\{ e^{2 \pi i n \cdot t}: n \in \bb
Z^N \},$ defines the Hilbert space isomorphism. This allows the
identification, $B(\ell^2(\bb Z^N)) = B(L^2([0,1]^N))$ and with
respect to this identification, $VN(\bb Z^N) = \{ M_f : f \in
L^{\infty}([0,1]^N) \},$ where $M_f$ denotes the operator of
multiplication by the function $f.$  The operator $M_f$ is often
called the {\em Laurent operator with symbol f,} especially in the
case $N=1.$  Thus, the above theorem, specialized to the case of $\bb
Z^N,$ shows that if a Laurent operator can be paved, then it can be
paved by syndetic sets in the sense of the above theorem.

Thus to find a Laurent operator that is not pavable and hence show
that the Kadison-Singer problem has a negative answer, it would be enough to
find a Laurent operator that can not be syndetically paved in the
above sense.
\end{remark}

We now turn our attention to the case of Bessel sequences in Hilbert spaces.
We need one observation.

\begin{lemma} Let $\{ f_i \}_{i \in I}$ be a Bessel sequence in a Hilbert space and let $A \subseteq I.$  Then $\{ f_i \}_{i \in A}$ is a Riesz basic sequence if and only if there exists $\delta > 0$ such that the operator inequality $\delta P_A \le P_A( \langle f_j,f_i \rangle ) P_A$ holds on $B(\ell^2(I)).$
\end{lemma}

We can now prove the results announced in the beginning.

{\em Proof of Theorem~\ref{thm1}.} Let $X = (x_{g,h})= ( \langle f_h,
f_g \rangle)$ and note that since $\{ f_g \}_{g \in G}$ is a Bessel
sequence, $X$ is the matrix of a bounded, positive operator on
$B(\ell^2(G)).$  The condition $\langle f_{g_1}, f_{g_2} \rangle =
\langle f_{h_1}, f_{h_2} \rangle$ whenever $g_1g_2^{-1} =
h_1h_2^{-1},$ guarantees that $X \in VN(G).$ To see this note that $X
\in VN(G)$ if and only if for every $g \in G,$ $E(X\lambda(g))$ is a
scalar multiple of the identity. We have that 
\[ \langle X\lambda(g)e_h,e_h \rangle = x_{h,gh} = \langle f_{gh},
f_{h} \rangle, \]
but for any two values of $h,$ we have that $(gh_1)h_1^{-1} =
(gh_2)h_2^{-1}$ and hence, $x_{h_1,gh_1} = x_{h_2,gh_2}.$ Thus, the
diagonal matrix $E(X\lambda(g))$ is some constant $c_g$ and $X \sim
\sum_{g \in G} c_g \lambda(g).$

Now since $\{ f_g: g \in A_l \}$ is a Bessel sequence for $1 \le l \le L,$ there exists $\delta > 0,$ such that $\delta P_{A_l} \le P_{A_l}XP_{A_l},$ for all $l.$
Fix any $p \in K(\beta G)$ and apply $\psi_p$ to this operator inequality, to obtain $\delta P_{S_l} = \delta \psi_p(P_{A_l}) \le \psi_p(P_{A_l}XP_{A_l}) = P_{S_l}XP_{S_l},$ where each $S_l$ is either empty or syndetic.

The inequality $\delta P_{S_l} \le P_{S_l}XP_{S_l}$ guarantees that $\{f_g: g \in S_l \}$ is a Riesz basic sequence and so the proof is complete.

{\em Proof of Theorem~\ref{thm2}.} We shall apply Theorem~\ref{thm1} with $G= \bb Z^N.$ Note that for any $n,m \in \bb Z^N,$ we have that
\[ \langle f_n,f_m \rangle = \int_A e_{n-m}(t) dt= \langle f_{n-m}, f_0 \rangle, \]
where $dt$ denotes the $N$-dimensional Lebesgue measure. From this identity it readily follows that the Bessel sequence $\{ f_n: n \in \bb Z^N \}$ satisfies the necessary group invariance condition. So the hypotheses of Theorem~\ref{thm1} are met.

\begin{remark} Let $S \subseteq \bb Z^N$ be a syndetic set and let $F
  \subseteq \bb Z^N$ be a finite set such that $F + S = \bb Z^N.$ For
  $n=(n_1,...,n_N) \in \bb Z^N$ set $\|n\| = \max \{ |n_1|,..., |n_N|
  \},$ i.e., the supremum norm. Let $M = \max \{ \|m\|: m \in F \}+1,$
  then given any $n \in S,$ there must exist $m \in S, m \ne n,$ with
  $\|n-m\| \le M.$ This in turn implies that $S$ will have positive
  ``density'' for most reasonable definitions of density. In particular, for
  $N=1$ every syndetic set has positive lower Beurling density.
Combining these observations with our earlier remark, we see that if
a Laurent operator can be paved, then it can be paved with syndetic
sets and consequently with sets of positive density.
\end{remark}

\begin{remark}  In the case $G= \bb Z,$ Halpern, Kaftal and
  Weiss\cite{HKW} studied {\em uniform paving} of Laurent
  operators. That is, they asked when $M_f \in B(\ell^2(\bb Z))$ with
  $E(M_f) =0$ had
  the property that for each $\epsilon >0,$ there existed an $M>0,$
  such that if one let $A = \{ Mn: n \in \bb Z \},$ then
  $\|P_AM_fP_A\| < \epsilon \|M_f\|.$  They proved that this occurs if
  and only if the symbol $f$ could be chosen Riemann integrable.

In contrast our result shows that if paving of a Laurent operator can
be done at all, then the set $A$ can be chosen syndetic and hence,
while not uniform, it will have a ``bounded gap,'' i.e., it must
select at least one element from each interval of length $M,$ for some $M.$
\end{remark}

This last remark lends credence to the belief that the symbols for which $M_f$ is
pavable are perhaps a proper subspace of $L^{\infty}([0,1]).$ Perhaps,
the symbols that can be paved are obtained from Riemann integrable
functions via some randomization method.

\end{document}